\documentclass{gtart_h}

\def\ifplaintex{\expandafter\ifx\csname documentclass\endcsname\relax}

\def\gtp{{\mathsurround=0pt\it $\cal G\mskip-2mu$eometry \&\ 
$\cal T\!\!$opology $\cal P\!$ublications}}  

\def\recd{{\small Received:\qua\receiveddate\ifx\reviseddate\relax
\else\qquad Revised:\qua\reviseddate\fi\par}} 

\def\lognumber#1{\def\thelognumber{#1}}
\def\volumenumber#1{\def\thevolumenumber{#1}}
\def\volumeyear#1{\def\thevolumeyear{#1}}
\def\papernumber#1{\def\thepapernumber{#1}}
\def\pagenumbers#1#2{\def\startpage{#1}\def\finishpage{#2}}
\def\published#1{\def\publishdate{#1}}

\def\received#1{\def\receiveddate{#1}}
\def\revised#1{\def\reviseddate{#1}}
\def\accepted#1{\def\accepteddate{#1}}

\long\def\asciiabstract#1{\long\def\theasciiabstract{#1}}
\def\asciikeywords#1{\def\theasciikeywords{#1}}

\let\\\par\let\thelognumber\relax\let\thevolumenumber\relax
\let\thepapernumber\relax\let\thevolumeyear\relax\let\startpage\relax
\let\finishpage\relax\let\publishdate\relax\let\receiveddate\relax
\let\reviseddate\relax\let\accepteddate\relax\let\theasciititle\relax
\let\theasciiauthors\relax
\let\theasciiabstract\relax\let\theasciikeywords\relax

\let\theasciiemail\relax

\ifplaintex
\font\logobig=cmssbx10 scaled 3836
\font\logomed=cmssbx10 scaled 2557
\else
\font\logobig=cmssbx10 scaled 4200
\font\logomed=cmssbx10 scaled 2800
\fi

\long\def\makeagttitle{   
\count0=\startpage
\agt\hfill      
\hbox to 45truept{\vbox to 0pt{\vglue -13truept{\logomed A\kern -.37em{\logobig 
T}\kern -.38em G}\vss}\hss}
\break
{\small Volume \thevolumenumber\ (\thevolumeyear)
\startpage--\finishpage\nl
Published: \publishdate}

\vglue .25truein

{\parskip=0pt\leftskip 0pt plus
1fil\def\\{\par\smallskip}{\Large\bf\thetitle}\par\medskip} \vglue
0.05truein

{\parskip=0pt\leftskip 0pt plus 1fil\def\\{\par}{\sc\theauthors}
\par\medskip}%
 
\vglue 0.03truein 

{\small\leftskip 25truept\rightskip 25truept{\bf Abstract}\stdspace\theabstract

{\bf AMS Classification}\stdspace\theprimaryclass
\ifx\thesecondaryclass\relax\else; \thesecondaryclass\fi\par
{\bf Keywords}\stdspace \thekeywords\par}\vglue 7truept

}   

\ifplaintex
\font\phead=cmsl9 scaled 950
\font\pnum=cmbx10 scaled 913
\font\pfoot=cmsl9 scaled 950
\headline{\vbox to 0pt{\vskip -4.5mm\line{\small\phead\ifnum
\count0=\startpage ISSN 1472-2739 (on-line) 1472-2747 (printed)
\hfill {\pnum\folio}\else\ifodd\count0\def\\{ }%
\ifx\theshorttitle\relax\thetitle\else\theshorttitle\fi\hfill{\pnum\folio}
\else\def\\{ and }{\pnum\folio}\hfill\ifx\theshortauthors\relax\theauthors
\else\theshortauthors\fi\fi\fi}\vss}}
\footline{\vbox to 0pt{\vglue 0mm\line{\small\pfoot\ifnum\count0=\startpage
\copyright\ \gtp\hfill\else
\agt, Volume \thevolumenumber\ (\thevolumeyear)\hfill\fi}\vss}}
\else
\font=cmsl9 scaled 1050
\font\lnum=cmbx10 
\font=cmsl9 scaled 1050
\makeatletter
\def\@oddhead{{\small\lhead\ifnum\count0=\startpage ISSN 1472-2739 
(on-line) 1472-2747 (printed)\hfill {\lnum\number\count0}\else\ifodd\count0
\def\\{ }\ifx\theshorttitle\relax \thetitle \else\theshorttitle\fi\hfill
{\lnum\number\count0}\else\def\\{ and }{\lnum\number\count0}
\hfill\ifx\theshortauthors\relax 
\theauthors\else\theshortauthors\fi\fi\fi}}\def\@evenhead{\@oddhead}
\def\@oddfoot{\small\lfoot\ifnum\count0=\startpage\copyright\ \gtp\hfill\else
\agt, Volume \thevolumenumber\ (\thevolumeyear)\hfill\fi}
\def\@evenfoot{\@oddfoot}
\makeatother
\fi
\let\maketitlepage\makeagttitle
\let\makeshorttitle\maketitlepage
\let\maketitle\maketitlepage

\newwrite\gtoutfile
\long\gdef\makeheadfile{  
{\def\\{, }\def\s{ }
\immediate\openout\gtoutfile head.xxx
\immediate\write\gtoutfile{Proxy-for: \ifx\theasciiauthors\relax
\theauthors\else\theasciiauthors\fi\s<\ifx\theasciiemail\relax\theemail\else\theasciiemail\fi>}
\immediate\write\gtoutfile{\noexpand\\}
\immediate\write\gtoutfile{Authors: \ifx\theasciiauthors\relax
\theauthors\else\theasciiauthors\fi}
{\def\\{ }\immediate\write\gtoutfile{Title: \ifx\theasciititle\relax
\thetitle\else\theasciititle\fi}}
\immediate\write\gtoutfile{Subj-class: GT or SG, GR etc}
\immediate\write\gtoutfile{MSC-class: \theprimaryclass\ifx\thesecondaryclass\relax\else, \thesecondaryclass\fi}
\immediate\write\gtoutfile{Journal-ref: Algebraic and Geometric Topology \thevolumenumber\s
(\thevolumeyear) \startpage-\finishpage}
\immediate\write\gtoutfile{Comments: Published by Algebraic and
Geometric Topology at}
\immediate\write\gtoutfile{\s\s\s  http://www.maths.warwick.ac.uk/agt/AGTVol\thevolumenumber/agt-\thevolumenumber-\thepapernumber.abs.html}
\immediate\write\gtoutfile{\noexpand\\}
\immediate\write\gtoutfile{}
\ifx\theasciiabstract\relax
\immediate\write\gtoutfile{\theabstract}\else
\immediate\write\gtoutfile{\theasciiabstract}\fi
\immediate\write\gtoutfile{}
\immediate\write\gtoutfile{\noexpand\\}
\immediate\write\gtoutfile{}
\immediate\closeout\gtoutfile}}  

\def\maketitlepage{\makeagttitle\makeheadfile}
\let\makeshorttitle\maketitlepage
\let\maketitle\maketitlepage

\lognumber{7}
\volumenumber{4}
\volumeyear{2004}
\papernumber{7}
\published{6 March 2004}
\pagenumbers{95}{119}
\received{29 September 2003}
\revised{1 March 2004}
\accepted{4 March 2004}

\usepackage{rotating}
\usepackage{amsmath,amssymb}
\usepackage[all]{xy}
\xyoption{v2}

\newtheorem{theorem}{Theorem}
\theoremstyle{definition}
\newtheorem{definition}{Definition}

\newcommand{\bcal}[1]{\mbox{\boldmath${\cal {#1}}$}}
\font\smallboldcal=cmbsy7 scaled 1200
\newcommand{\bbcal}[1]{\mbox{$\mbox{\smallboldcal #1}$}}

\newenvironment{fake}{\relax}{\relax}

\begin{document}

\title{Enrichment over iterated monoidal categories}
\author{Stefan Forcey}
\address{Department of Mathematics, Virginia Tech\\460 McBryde Hall, 
Blacksburg, VA 24060, USA}
\email{sforcey@math.vt.edu}
\asciikeywords{Loop spaces, enriched categories, n-categories, 
iterated monoidal categories}

\begin{abstract}
Joyal and Street note in their paper on braided monoidal categories
\cite{JS} that the 2--category ${\cal V}$--Cat of categories enriched
over a braided monoidal category ${\cal V}$ is not itself braided in
any way that is based upon the braiding of ${\cal V}$.  The exception
that they mention is the case in which ${\cal V}$ is symmetric, which
leads to ${\cal V}$--Cat being symmetric as well.  The symmetry in
${\cal V}$--Cat is based upon the symmetry of ${\cal V}$.  The
motivation behind this paper is in part to describe how these facts
relating ${\cal V}$ and ${\cal V}$--Cat are in turn related to a
categorical analogue of topological delooping.  To do so I need to
pass to a more general setting than braided and symmetric categories
--- in fact the $k$--fold monoidal categories of Balteanu et al in
\cite{Balt}. It seems that the analogy of loop spaces is a good guide
for how to define the concept of enrichment over various types of
monoidal objects, including $k$--fold monoidal categories and their
higher dimensional counterparts.  The main result is that for ${\cal
V}$ a $k$--fold monoidal category, ${\cal V}$--Cat becomes a
$(k-1)$--fold monoidal $2$--category in a canonical way.  In the next
paper I indicate how this process may be iterated by enriching over
${\cal V}$--Cat, along the way defining the 3--category of categories
enriched over ${\cal V}$--Cat.  In future work I plan to make precise
the $n$--dimensional case and to show how the group completion of the
nerve of ${\cal V}$ is related to the loop space of the group
completion of the nerve of ${\cal V}$--Cat.

This paper is an abridged  version of \cite{forcey1}.
\end{abstract}

\asciiabstract{Joyal and Street note in their paper on braided
monoidal categories [Braided tensor categories, Advances in
Math. 102(1993) 20-78] that the 2-category V-Cat of categories
enriched over a braided monoidal category V is not itself braided in
any way that is based upon the braiding of V.  The exception that they
mention is the case in which V is symmetric, which leads to V-Cat
being symmetric as well.  The symmetry in V-Cat is based upon the
symmetry of V.  The motivation behind this paper is in part to
describe how these facts relating V and V-Cat are in turn related to a
categorical analogue of topological delooping.  To do so I need to
pass to a more general setting than braided and symmetric categories
-- in fact the k-fold monoidal categories of Balteanu et al in
[Iterated Monoidal Categories, Adv. Math. 176(2003) 277-349]. It seems
that the analogy of loop spaces is a good guide for how to define the
concept of enrichment over various types of monoidal objects,
including k-fold monoidal categories and their higher dimensional
counterparts.  The main result is that for V a k-fold monoidal
category, V-Cat becomes a (k-1)-fold monoidal 2-category in a
canonical way.  In the next paper I indicate how this process may be
iterated by enriching over V-Cat, along the way defining the
3-category of categories enriched over V-Cat.  In future work I plan
to make precise the n-dimensional case and to show how the group
completion of the nerve of V is related to the loop space of the group
completion of the nerve of V-Cat.  

This paper is an abridged version of `Enrichment as categorical
delooping I: Enrichment over iterated monoidal categories',
math.CT/0304026}

\begin{fake}\end{fake}
\maketitle

\SloppyCurves
\section{Introduction}

    A major goal of higher dimensional category theory is to discover
    ways of exploiting the connections between homotopy coherence and
    categorical coherence. Stasheff \cite{Sta} and Mac Lane \cite{Mac}
    showed that monoidal categories are precisely analogous to 1--fold
    loop spaces. There is a similar connection between symmetric
    monoidal categories and infinite loop spaces. The first step in
    filling in the gap between 1 and infinity was made in \cite{ZF}
    where it is shown that the group completion of the nerve of a
    braided monoidal category is a 2--fold loop space.  In \cite{Balt}
    the authors finished this process by, in their words, ``pursuing
    an analogy to the tautology that an $n$--fold loop space is a loop
    space in the category of $(n-1)$--fold loop spaces.'' The first
    thing they focus on is the fact that a braided category is a
    special case of a carefully defined 2--fold monoidal
    category. Based on their observation of the correspondence between
    loop spaces and monoidal categories, they iteratively define the
    notion of $n$--fold monoidal category as a monoid in the category
    of $(n-1)$--fold monoidal categories.  In their view ``monoidal''
    functors should be defined in a more ``lax'' way than is usual in
    order to avoid strict commutativity of 2--fold and higher monoidal
    categories. In \cite{Balt} a symmetric category is seen as a
    category that is $n$--fold monoidal for all $n$.
    
    The main result in \cite{Balt} is that their definition of
    iterated monoidal categories exactly corresponds to $n$--fold loop
    spaces for all $n$. They show that the group completion of the
    nerve of an $n$--fold monoidal category is an $n$--fold loop
    space. Then they describe an operad in the category of small
    categories which parameterizes the algebraic structure of an
    iterated monoidal category.  They show that the nerve of this
    categorical operad is a topological operad which is equivalent to
    the little $n$--cubes operad. This latter operad, as shown in
    \cite{BV1} and \cite{May}, characterizes the notion of $n$--fold
    loop space.  Thus the main result in \cite{Balt} is a categorical
    characterization of $n$--fold loop spaces.
    
    The present paper pursues the hints of a categorical delooping
    that are suggested by the facts that for a symmetric category, the
    2--category of categories enriched over it is again symmetric,
    while for a braided category the 2--category of categories
    enriched over it is merely monoidal. Section 2 reviews enrichment.
    Section 3 goes over the recursive definition of the $k$--fold
    monoidal categories of \cite{Balt}, altered here to include a
    coherent associator.  The immediate question is whether the
    delooping phenomenon happens in general for these $k$--fold
    monoidal categories.  The answer is yes, once enriching over a
    $k$--fold monoidal category is carefully defined in Section 4,
    where we see that all the information included in the axioms for
    the $k$--fold category is exhausted in the process. The definition
    also provides for iterated delooping as is previewed in Section 5.
    
    I have organized the paper so that sections can largely stand
    alone, so please skip them when able, and forgive redundancy when
    it occurs. Thanks to my advisor, Frank Quinn, for inspirational
    suggestions.  Thanks to {\Xy-pic} for the diagrams. Thanks
    especially to the authors of \cite{Balt} for making their source
    available--I learned and borrowed from their use of {\LaTeX} as
    well as from their insights into the subject matter.

\section{Review of categories enriched over a monoidal category}

   In this section I briefly review the definition of a category
  enriched over a monoidal category ${\cal V}$.  Enriched functors and
  enriched natural transformations make the collection of enriched
  categories into a 2-category ${\cal V}$-Cat.  This section is not
  meant to be complete. It is included due to, and its contents
  determined by, how often the definitions herein are referred to and
  followed as models in the rest of the paper. The definitions and
  proofs can be found in more or less detail in \cite{Kelly} and
  \cite{EK1} and of course in \cite{MacLane}.
    
    \begin{definition} For our purposes a {\it monoidal category} is a
      category ${\cal V}$ together with a functor $\otimes: {\cal
      V}\times{\cal V}\to{\cal V}$ and an object $I$ such that:
      \begin{enumerate} 

\item $\otimes$ is associative up to the coherent natural
      transformations $\alpha$. The coherence axiom is given by the
      commuting pentagon: 
\vspace{2mm}
\begin{center} \resizebox{4.5in}{!}{
              $$
              \xymatrix@C=-35pt{
              &((U\otimes V)\otimes W)\otimes X \text{ }\text{ }
              \ar[rr]^{ \alpha_{UVW}\otimes 1_{X}}
              \ar[dl]^{ \alpha_{(U\otimes V)WX}}
              &\text{ }\text{ }\text{ }\text{ }\text{ }\text{ }\text{ }\text{ }\text{ }&\text{ }\text{ }(U\otimes (V\otimes W))\otimes X
              \ar[dr]^{ \alpha_{U(V\otimes W)X}}&\\
              (U\otimes V)\otimes (W\otimes X)
              \ar[drr]|{ \alpha_{UV(W\otimes X)}}
              &&&&U\otimes ((V\otimes W)\otimes X)
              \ar[dll]|{ 1_{U}\otimes \alpha_{VWX}}
              \\&&U\otimes (V\otimes (W\otimes X))&&&
             }
             $$
             }
	                    \end{center}
	                    
    \item $I$ is a strict $2$-sided unit for $\otimes$.
    \end{enumerate}
    \end{definition}
    \begin {definition} \label{V:Cat} A (small) ${\cal V}$ {\it--Category} ${\cal A}$ is a set $\left|{\cal A}\right|$ of 
    {\it objects}, 
    a {\it hom-object} ${\cal A}(A,B) \in \left|{\cal V}\right|$ for
    each pair of objects of ${\cal A}$, a family of {\it composition morphisms} $M_{ABC}:{\cal A}(B,C)
    \otimes{\cal A}(A,B)\to{\cal A}(A,C)$ for each triple of objects, and an {\it identity element} $j_{A}:I\to{\cal A}(A,A)$ for each object.
    The composition morphisms are subject to the associativity axiom which states that the following pentagon commutes
\vspace{-3mm}
          	          \begin{center}
    	          \resizebox{4.5in}{!}{
          $$
          \xymatrix@C=-22pt{
          &({\cal A}(C,D)\otimes {\cal A}(B,C))\otimes {\cal A}(A,B)\text{ }\text{ }
          \ar[rr]^{\scriptstyle \alpha}
          \ar[dl]^{\scriptstyle M \otimes 1}
          &\text{ }\text{ }\text{ }\text{ }\text{ }\text{ }\text{ }\text{ }\text{ }\text{ }\text{ }\text{ }\text{ }\text{ }\text{ }&\text{ }\text{ }{\cal A}(C,D)\otimes ({\cal A}(B,C)\otimes {\cal A}(A,B))
          \ar[dr]^{\scriptstyle 1 \otimes M}&\\
          {\cal A}(B,D)\otimes {\cal A}(A,B)
          \ar[drr]^{\scriptstyle M}
          &&&&{\cal A}(C,D)\otimes {\cal A}(A,C)
          \ar[dll]^{\scriptstyle M}
          \\&&{\cal A}(A,D))&&&
          }$$
          }
         	                    \end{center}
and to the unit axioms which state that both the triangles in the following diagram commute.
     $$
      \xymatrix{
      I\otimes {\cal A}(A,B)
      \ar[rrd]^{=}
      \ar[dd]_{j_{B}\otimes 1}
      &&&&{\cal A}(A,B)\otimes I 
      \ar[dd]^{1\otimes j_{A}}
      \ar[lld]^{=}\\
      &&{\cal A}(A,B)\\
      {\cal A}(B,B)\otimes {\cal A}(A,B)
      \ar[rru]^{M_{ABB}}
      &&&&{\cal A}(A,B)\otimes {\cal A}(A,A)
      \ar[llu]^{M_{AAB}}
      }
   $$
   \end{definition}

    In general a ${\cal V}$--category is directly analogous to an (ordinary) category enriched over $\mathbf{Set}.$  
    If ${\cal V} = \mathbf{Set}$ then these diagrams are the usual category axioms.

   \begin{definition} \label{enriched:funct} For ${\cal V}$--categories 
   ${\cal A}$ and ${\cal B}$, a ${\cal V}$--$functor$ $T:{\cal A}\to{\cal B}$ is a function
    $T:\left| {\cal A} \right| \to \left| {\cal B} \right|$ and a family of 
    morphisms $T_{AB}:{\cal A}(A,B) \to {\cal B}(TA,TB)$ in ${\cal V}$ indexed by 
    pairs $A,B \in \left| {\cal A} \right|$.
    The usual rules for a functor that state $T(f \circ g) = Tf \circ Tg$ 
    and $T1_{A} = 1_{TA}$ become in the enriched setting, respectively, the commuting diagrams
   $$
    \xymatrix{
    &{\cal A}(B,C)\otimes {\cal A}(A,B)
    \ar[rr]^{\scriptstyle M}
    \ar[d]^{\scriptstyle T \otimes T}
    &&{\cal A}(A,C)
    \ar[d]^{\scriptstyle T}&\\
    &{\cal B}(TB,TC)\otimes {\cal B}(TA,TB)
    \ar[rr]^{\scriptstyle M}
    &&{\cal B}(TA,TC)
    }
   $$
  and
   $$
    \xymatrix{
    &&{\cal A}(A,A)
    \ar[dd]^{\scriptstyle T_{AA}}\\
    I
    \ar[rru]^{\scriptstyle j_{A}}
    \ar[rrd]_{\scriptstyle j_{TA}}\\
    &&{\cal B}(TA,TA).
    }
   $$
  ${\cal V}$--functors can be composed to form a category called ${\cal V}$--Cat. This category
  is actually enriched over $\mathbf{Cat}$, the category of (small) categories with cartesian product. 
     \end{definition}

  \begin{definition} \label{enr:nat:trans}
  For ${\cal V}$--functors $T,S:{\cal A}\to{\cal B}$ a ${\cal V}$--{\it natural  
   transformation} $\alpha:T \to S:{\cal A} \to {\cal B}$
  is an $\left| {\cal A} \right|$--indexed family of 
  morphisms $\alpha_{A}:I \to {\cal B}(TA,SA)$ satisfying the ${\cal V}$--naturality
  condition expressed by the commutativity of the following hexagonal diagram:
  	            	          \begin{center}
      	          \resizebox{4.5in}{!}{
   $$
    \xymatrix{
    &I \otimes {\cal A}(A,B)
    \ar[rr]^-{\scriptstyle \alpha_{B} \otimes T_{AB}}
    &&{\cal B}(TB,SB) \otimes {\cal B}(TA,TB)
    \ar[rd]^-{\scriptstyle M}
  \\
    {\cal A}(A,B)
    \ar[ru]^{=}
  \ar[rd]_{=}
    &&&&{\cal B}(TA,SB)
  \\
    &{\cal A}(A,B) \otimes I
    \ar[rr]_-{\scriptstyle S_{AB} \otimes \alpha_{A}}
    &&{\cal B}(SA,SB) \otimes {\cal B}(TA,SA)
    \ar[ru]^-{\scriptstyle M}
    }
   $$
  }
         	                    \end{center}
   \end{definition}
   
   For two ${\cal V}$--functors  $T,S$ to be equal is to say $TA = SA$ for all $A$ 
   and for the ${\cal V}$--natural isomorphism $\alpha$ between them to have 
   components $\alpha_{A} = j_{TA}$. This latter implies equality of the hom--object morphisms: 
   $T_{AB} = S_{AB}$ for all pairs of objects. The implication is seen by combining the second diagram in 
   Definition \ref{V:Cat} with all the diagrams
  in Definitions \ref{enriched:funct} and \ref{enr:nat:trans}.
   
  We want to check that ${\cal V}$--natural transformations can be composed so that  ${\cal V}$--categories, ${\cal V}$--functors
  and ${\cal V}$--natural transformations form a 2--category. First the vertical composite of ${\cal V}$--natural transformations corresponding to the 
  picture
    $$
  \xymatrix@R-=16pt{
  &\ar@{=>}[d]^{\alpha}\\
  {\cal A}
  \ar@/^2pc/[rr]^T
  \ar[rr]_>>>>>S
  \ar@/_2pc/[rr]_R
  &\ar@{=>}[d]^{\beta}
  &{\cal B}\\
  &\\
  }
  $$
  has components given by 
  $(\beta \circ \alpha)_{A} = \xymatrix{
                                  I \cong I \otimes I 
                                  \ar[d]_{\beta_{A} \otimes \alpha_{A}}
                                  \\{\cal B}(SA,RA) \otimes {\cal B}(TA,SA)
                                  \ar[d]_M
                                  \\{\cal B}(TA,RA).}$
  
  The reader should check that this composition produces a valid
  ${\cal V}$--natural transformation. Associativity of composition
  also follows from the pentagonal axioms.  The identity 2-cells are
  the identity ${\cal V}$-natural transformations ${\textbf 1}_{Q} : Q
  \to Q : {\cal B} \to {\cal C}.$ These are formed from the unit
  morphisms in ${\cal V}$: ${({\textbf 1}_{Q})}_{B} = j_{QB}$.

  In order to define composition of all allowable pasting diagrams in
  the 2-category, we need to define the composition described by left
  and right whiskering diagrams. The first picture shows a 1-cell
  (${\cal V}$--functor) following a 2-cell (${\cal V}$--natural
  transformation).  These are composed to form a new 2-cell as follows
    $$
  \xymatrix@R-=3pt{
  &\ar@{=>}[dd]^{\alpha}\\
  {\cal A}
  \ar@/^1pc/[rr]^T
  \ar@/_1pc/[rr]_S
  &&{\cal B}
  \ar[rr]^Q
  &&{\cal C}
  \\
  &\\
  } \text{   is composed to become } \xymatrix@R-=3pt{
  &\ar@{=>}[dd]^{Q\alpha}\\
  {\cal A}
  \ar@/^1pc/[rr]^{QT}
  \ar@/_1pc/[rr]_{QS}
  &&{\cal C}
  \\
  &\\
  }
  $$
    where $QT$ and $QS$ are given by the usual compositions of their
  set functions and morphisms in ${\cal V}$, and $Q\alpha$ has
  components given by $(Q\alpha)_{A} = \xymatrix{ I
  \ar[d]_{\alpha_{A}} \\{\cal B}(TA,SA) \ar[d]_{Q_{TA,SA}} \\{\cal
  C}(QTA,QSA).}$
                                  
  The second picture shows a 2-cell following a 1-cell. These are
  composed as follows
  $$
  \xymatrix@R-=3pt{
  &&&\ar@{=>}[dd]^{\alpha}\\
  {\cal D}
  \ar[rr]^P
  &&{\cal A}
  \ar@/^1pc/[rr]^T
  \ar@/_1pc/[rr]_S
  &&{\cal B}
  \\
  &&&\\
  } \text{   is composed to become } \xymatrix@R-=3pt{
  &\ar@{=>}[dd]^{\alpha P}\\
  {\cal D}
  \ar@/^1pc/[rr]^{TP}
  \ar@/_1pc/[rr]_{SP}
  &&{\cal B}
  \\
  &\\
  }
  $$
  where $\alpha P$ has components given by $(\alpha P)_{D} =
  \alpha_{PD}$.  What we have developed here are the partial functors
  of the composition morphism implicit in enriching over
  $\mathbf{Cat}$. The partial functors can be combined to make the
  functor of two variables as shown in \cite{EK1}.

  Having ascertained that we have a 2--category we review the
   morphisms between two such things.  A {\it 2--functor} $F:U\to V$
   sends objects to objects, 1--cells to 1--cells, and 2--cells to
   2--cells and preserves all the categorical structure. A {\it
   2--natural transformation} $\theta:F\to G:U\to V$ is a function
   that sends each object $A \in U$ to a 1--cell $\theta_{A}:FA\to GA$
   in $V$ in such a way that for each 2--cell in $U$ the compositions
   of the following diagrams are equal in $V$.
   $$
   \xymatrix@R-=3pt{
   &\ar@{=>}[dd]^{F\alpha}\\
    FA
   \ar@/^1pc/[rr]^{Ff}
   \ar@/_1pc/[rr]_{Fg}
   &&FB
   \ar[rr]^{\theta_B}
   &&GB
   \\
  &\\
  }
  =
  \xymatrix@R-=3pt{
  &&&\ar@{=>}[dd]^{G\alpha}\\
  FA
  \ar[rr]^{\theta_A}
  &&GA
  \ar@/^1pc/[rr]^{Gf}
  \ar@/_1pc/[rr]_{Gg}
  &&GB
  \\
  &&&\\
  }
  $$

\section{$k$-fold monoidal categories}
    
    In this section I closely follow the authors of \cite{Balt} in defining a notion of iterated monoidal category.
    For those readers familiar with that source, note that I vary from their definition only 
    by including associativity 
    up to coherent natural isomorphisms.
    This includes changing the basic picture from monoids to something that is a monoid only
    up to a monoidal natural transformation.
    We (and in this section ``we'' is not merely imperial, since so much is directly from \cite{Balt})
    start by
    defining a slightly nonstandard variant
    of monoidal functor. It is usually required in a definition of monoidal functor
    that $\eta$ be an isomorphism. The authors of \cite {Balt} note that it is crucial not to make this
    requirement. 
    
    \begin{definition} A {\it monoidal functor} $(F,\eta) :{\cal
    C}\to{\cal D}$ between monoidal categories consists of a functor
    $F$ such that $F(I)=I$ together with a natural transformation
    $$
    \eta_{AB}:F(A)\otimes F(B)\to F(A\otimes B),
    $$
    which satisfies the following conditions.
    \begin{enumerate}
    \item Internal Associativity: The following diagram commutes. 
    $$
    \diagram
    (F(A)\otimes F(B))\otimes F(C)
    \rrto^{\eta_{AB}\otimes 1_{F(C)}}
    \dto^{\alpha}
    &&F(A\otimes B)\otimes F(C)
    \dto^{\eta_{(A\otimes B)C}}\\
    F(A)\otimes (F(B)\otimes F(C))
    \ar[d]^{1_{F(A)}\otimes \eta_{BC}}
    &&F((A\otimes B)\otimes C)
    \ar[d]^{F\alpha}\\
    F(A)\otimes F(B\otimes C)
    \rrto^{\eta_{A(B\otimes C)}}
    &&F(A\otimes (B\otimes C))
    \enddiagram
    $$
    \item Internal Unit Conditions: $\eta_{AI}=\eta_{IA}=1_{F(A)}.$
    \end{enumerate}
    \end{definition}
    Given two monoidal functors $(F,\eta) :{\cal C}\to{\cal D}$ and $(G,\zeta) 
    :{\cal D}\to{\cal E}$,
    we define their composite to be the monoidal functor $(GF,\xi) :
    {\cal C}\to{\cal E}$, where
    $\xi$ denotes the composite
    $$
    \diagram
    GF(A)\otimes GF(B)\rrto^{\zeta_{F(A)F(B)}} 
    && G\bigl(F(A)\otimes F(B)\bigr)\rrto^{G(\eta_{AB})}
    &&GF(A\otimes B).
    \enddiagram
    $$
    It is easy to verify that $\xi$ satisfies the internal associativity condition above by subdividing the 
    necessary commuting diagram into two regions that commute by the axioms for $\eta$ and $\zeta$ respectively 
    and two that commute due to their naturality.
     $\mathbf{MonCat}$ is the monoidal category of monoidal categories and monoidal
    functors, with the usual Cartesian product as in ${\mathbf{Cat}}$. 
    
    A {\it monoidal natural transformation} $\theta:(F, \eta) \to (G, \zeta):{\cal D}\to{\cal E}$  is a 
    natural transformation $\theta: F\to G$ between the underlying ordinary functors of $F$ and $G$ such that the
    following diagram commutes
    $$
    \xymatrix{
    F(A)\otimes F(B)
    \ar[r]^{\eta}
    \ar[d]^{\theta_A \otimes \theta_B}
    &F(A\otimes B)
    \ar[d]^{\theta_{A \otimes B}}
    \\
    G(A)\otimes G(B)
    \ar[r]^{\zeta}
    &G(A\otimes B)
    }
    $$

    \begin{definition} For our purposes a $2${\it -fold monoidal category} is a tensor object, or pseudomonoid, 
    in $\mathbf{MonCat}$.
    This means that we are given a monoidal category $({\cal V},\otimes_1,\alpha^1,I)$ and a 
    monoidal functor
    $(\otimes_2,\eta):{\cal V}\times{\cal V}\to{\cal V}$ which satisfies:
    \begin{enumerate}
    \item External Associativity: the following diagram describes a monoidal natural transformation
    $\alpha^2$ in $\mathbf{MonCat}.$ 
    $$
    \xymatrix{
    {\cal V}\times{\cal V}\times{\cal V}
    \rrto^{(\otimes_2,\eta)\times 1_{\cal V}}
    \dto_(0.4){1_{\cal V}\times(\otimes_2,\eta)}
    && {\cal V}\times{\cal V}
    \dto^{(\otimes_2,\eta)}
    \ar@{=>}[dll]^{\alpha^2}
    \\
    {\cal V}\times{\cal V} 
    \rrto_{(\otimes_2,\eta)}
    &&
    {\cal V}
    }
    $$
    
    \item External Unit Conditions: the following diagram commutes in 
    $\mathbf{MonCat}$. 
    $$
    \diagram
    {\cal V}\times I 
    \rto^{\subseteq}
    \ddrto^{\cong}
    & {\cal V}\times{\cal V}
    \ddto^{(\otimes_2,\eta)}
    & I\times{\cal V}
    \lto_{\supseteq}
    \ddlto^{\cong}\\\\
    &{\cal V} 
    \enddiagram
    $$
    \item Coherence: The underlying natural transformation $\alpha^2$ satisfies
    the usual coherence pentagon.
       
    \end{enumerate}
    \end{definition}
    
    Explicitly this means that we are given a second associative binary operation
    $\otimes_2:{\cal V}\times{\cal V}\to{\cal V}$, for which $I$ is also a two-sided unit.
   We are also given a natural transformation
    $$
    \eta_{ABCD}: (A\otimes_2 B)\otimes_1 (C\otimes_2 D)\to
    (A\otimes_1 C)\otimes_2(B\otimes_1 D).
    $$
    The internal unit conditions give $\eta_{ABII}=\eta_{IIAB}=1_{A\otimes_2 B}$,
    while the external unit conditions give $\eta_{AIBI}=\eta_{IAIB}=1_{A\otimes_1 B}$.
    The internal associativity condition gives the commutative diagram:    
    \noindent
    	            	          \begin{center}
        	          \resizebox{5in}{!}{
                 $$
    \diagram
    ((U\otimes_2 V)\otimes_1 (W\otimes_2 X))\otimes_1 (Y\otimes_2 Z)
    \xto[rrr]^{\eta_{UVWX}\otimes_1 1_{Y\otimes_2 Z}}
    \ar[d]^{\alpha^1}
    &&&\bigl((U\otimes_1 W)\otimes_2(V\otimes_1 X)\bigr)\otimes_1 (Y\otimes_2 Z)
    \dto^{\eta_{(U\otimes_1 W)(V\otimes_1 X)YZ}}\\
    (U\otimes_2 V)\otimes_1 ((W\otimes_2 X)\otimes_1 (Y\otimes_2 Z))
    \dto^{1_{U\otimes_2 V}\otimes_1 \eta_{WXYZ}}
    &&&((U\otimes_1 W)\otimes_1 Y)\otimes_2((V\otimes_1 X)\otimes_1 Z)
    \ar[d]^{\alpha^1 \otimes_2 \alpha^1}
    \\
    (U\otimes_2 V)\otimes_1 \bigl((W\otimes_1 Y)\otimes_2(X\otimes_1 Z)\bigr)
    \xto[rrr]^{\eta_{UV(W\otimes_1 Y)(X\otimes_1 Z)}}
    &&& (U\otimes_1 (W\otimes_1 Y))\otimes_2(V\otimes_1 (X\otimes_1 Z))
    \enddiagram
    $$
                 }
         	                    \end{center}
    The external associativity condition gives the commutative diagram:
    \noindent
    	            	          \begin{center}
        	          \resizebox{5in}{!}{
 $$
    \diagram
    ((U\otimes_2 V)\otimes_2 W)\otimes_1 ((X\otimes_2 Y)\otimes_2 Z)
    \xto[rrr]^{\eta_{(U\otimes_2 V)W(X\otimes_2 Y)Z}}
    \ar[d]^{\alpha^2 \otimes_1 \alpha^2}
    &&& \bigl((U\otimes_2 V)\otimes_1 (X\otimes_2 Y)\bigr)\otimes_2(W\otimes_1 Z)
    \dto^{\eta_{UVXY}\otimes_2 1_{W\otimes_1 Z}}\\
    (U\otimes_2 (V\otimes_2 W))\otimes_1 (X\otimes_2 (Y\otimes_2 Z))
    \dto^{\eta_{U(V\otimes_2 W)X(Y\otimes_2 Z)}}
    &&&((U\otimes_1 X)\otimes_2(V\otimes_1 Y))\otimes_2(W\otimes_1 Z)
    \ar[d]^{\alpha^2}
    \\
    (U\otimes_1 X)\otimes_2\bigl((V\otimes_2 W)\otimes_1 (Y\otimes_2 Z)\bigr)
    \xto[rrr]^{1_{U\otimes_1 X}\otimes_2\eta_{VWYZ}}
    &&& (U\otimes_1 X)\otimes_2((V\otimes_1 Y)\otimes_2(W\otimes_1 Z))
    \enddiagram
    $$
    }
             	                    \end{center}
    The authors of \cite{Balt} remark that we have natural transformations
    $$
    \eta_{AIIB}:A\otimes_1 B\to A\otimes_2 B\qquad\mbox{ and }\qquad
    \eta_{IABI}:A\otimes_1 B\to B\otimes_2 A.
    $$
    If they had insisted a 2-fold monoidal category be a tensor object
    in the category of monoidal categories and {\it strictly
    monoidal\/} functors, this would be equivalent to requiring that
    $\eta=1$.  In view of the above, they note that this would imply
    $A\otimes_1 B = A\otimes_2 B = B\otimes_1 A$ and similarly for
    morphisms.
        
     Joyal and Street \cite{JS} considered a 
    similar concept to Balteanu, Fiedorowicz, Schw${\rm \ddot a}$nzl and Vogt's idea of 2--fold monoidal category.  
    The former pair required the natural transformation $\eta_{ABCD}$ 
    to be an isomorphism and showed that the resulting category is naturally 
    equivalent to a braided monoidal category. As explained in \cite{Balt}, given such a category one 
    obtains an equivalent braided monoidal category by discarding one of the two 
    operations, say $\otimes_2$, and defining the commutativity isomorphism for the 
    remaining operation $\otimes_1$ to be the composite
    $$
    \diagram
    A\otimes_1 B\rrto^{\eta_{IABI}} 
    && B\otimes_2 A\rrto^{\eta_{BIIA}^{-1}}
    && B\otimes_1 A.
    \enddiagram
    $$
    Just as in \cite{Balt} we now define a 2--fold monoidal functor
    between 2--fold monoidal categories $F:{\cal V}\to{\cal D}$. It is
    a functor together with two natural transformations:
    $$\lambda^1_{AB}:F(A)\otimes_1 F(B)\to F(A\otimes_1 B)$$
    $$\lambda^2_{AB}:F(A)\otimes_2 F(B)\to F(A\otimes_2 B)$$
    satisfying the same associativity and unit conditions as in the
    case of monoidal functors.  In addition we require that the
    following hexagonal interchange diagram commutes.   
   \noindent
    	            	          \begin{center}
        	          \resizebox{5in}{!}{
                 $$
    \diagram
    (F(A)\otimes_2 F(B))\otimes_1(F(C)\otimes_2 F(D))
    \xto[rrr]^{ {\eta_{F(A)F(B)F(C)F(D)}}}
    \dto^{ {\lambda^2_{AB}\otimes_1\lambda^2_{CD}}}
    &&&(F(A)\otimes_1 F(C))\otimes_2(F(B)\otimes_1 F(D))
    \dto^{ \lambda^1_{AC}\otimes_2\lambda^1_{BD}}\\
    F(A\otimes_2 B)\otimes_1 F(C\otimes_2 D)
    \dto^{ {\lambda^1_{(A\otimes_2 B)(C\otimes_2 D)}}}
    &&&F(A\otimes_1 C)\otimes_2 F(B\otimes_1 D)
    \dto^{ {\lambda^2_{(A\otimes_1 C)(B\otimes_1 D)}}}\\
    F((A\otimes_2 B)\otimes_1(C\otimes_2 D))
    \xto[rrr]^{ F(\eta_{ABCD})}
    &&& F((A\otimes_1 C)\otimes_2(B\otimes_1 D))
    \enddiagram
    $$
    }
             	                    \end{center}
    We can now define the category {\textbf{2}$\mathbf{ -MonCat}$ 
    of 2-fold monoidal categories and
    2-fold monoidal functors, and then define a 3-fold monoidal category
    as a tensor object in {\textbf{2}$\mathbf{ -MonCat}$.  From this
    point on, the iteration of this idea is straightforward and,
    paralleling the authors of \cite{Balt}, we arrive at the following
    definitions.
  
  \begin{definition} An $n${\it -fold monoidal category} is a category
          ${\cal V}$ with the following structure.  \begin{enumerate}
          \item There are $n$ distinct multiplications
          $$\otimes_1,\otimes_2,\dots, \otimes_n:{\cal V}
          \times{\cal V}\to{\cal V}$$
          for each of which the associativity pentagon commutes:
          \noindent\vspace{2mm}
  		          \begin{center}
  	          \resizebox{4.5in}{!}{
          $$
          \xymatrix@C=-35pt{
          &((U\otimes_i V)\otimes_i W)\otimes_i X \text{ }\text{ }
          \ar[rr]^{ \alpha^{i}_{UVW}\otimes_i 1_{X}}
          \ar[dl]^{ \alpha^{i}_{(U\otimes_i V)WX}}
          &\text{ }\text{ }\text{ }\text{ }\text{ }&\text{ }\text{ }(U\otimes_i (V\otimes_i W))\otimes_i X
          \ar[dr]^{ \alpha^{i}_{U(V\otimes_i W)X}}&\\
          (U\otimes_i V)\otimes_i (W\otimes_i X)
          \ar[drr]^{ \alpha^{i}_{UV(W\otimes_i X)}}
          &&&&U\otimes_i ((V\otimes_i W)\otimes_i X)
          \ar[dll]^{ 1_{U}\otimes_i \alpha^{i}_{VWX}}
          \\&&U\otimes_i (V\otimes_i (W\otimes_i X))&&&
          }
          $$
          }
  		                    \end{center}
          ${\cal V}$ has an object $I$ which is a strict unit for all
          the multiplications.  \item For each pair $(i,j)$ such that
          $1\le i<j\le n$ there is a natural transformation
          $$\eta^{ij}_{ABCD}: (A\otimes_j B)\otimes_i(C\otimes_j D)\to
          (A\otimes_i C)\otimes_j(B\otimes_i D).$$ \end{enumerate}
          These natural transformations $\eta^{ij}$ are subject to the
          following conditions: 

\begin{enumerate}  \item[(a)] Internal unit condition:
          $\eta^{ij}_{ABII}=\eta^{ij}_{IIAB}=1_{A\otimes_j B}$
          \item[(b)] External unit condition:
          $\eta^{ij}_{AIBI}=\eta^{ij}_{IAIB}=1_{A\otimes_i B}$
          \item[(c)] Internal associativity condition: The following
          diagram commutes.  \vspace{-2mm}\begin{center} \resizebox{5in}{!}{
	     $$
            \diagram ((U\otimes_j V)\otimes_i (W\otimes_j X))\otimes_i
            (Y\otimes_j Z) \xto[rrr]^{\eta^{ij}_{UVWX}\otimes_i
            1_{Y\otimes_j Z}} \ar[d]^{\alpha^i} &&&\bigl((U\otimes_i
            W)\otimes_j(V\otimes_i X)\bigr)\otimes_i (Y\otimes_j Z)
            \dto^{\eta^{ij}_{(U\otimes_i W)(V\otimes_i X)YZ}}\\
            (U\otimes_j V)\otimes_i ((W\otimes_j X)\otimes_i
            (Y\otimes_j Z)) \dto^{1_{U\otimes_j V}\otimes_i
            \eta^{ij}_{WXYZ}} &&&((U\otimes_i W)\otimes_i
            Y)\otimes_j((V\otimes_i X)\otimes_i Z) \ar[d]^{\alpha^i
            \otimes_j \alpha^i} \\ (U\otimes_j V)\otimes_i
            \bigl((W\otimes_i Y)\otimes_j(X\otimes_i Z)\bigr)
            \xto[rrr]^{\eta^{ij}_{UV(W\otimes_i Y)(X\otimes_i Z)}} &&&
            (U\otimes_i (W\otimes_i Y))\otimes_j(V\otimes_i
            (X\otimes_i Z)) \enddiagram
            $$
           }
	            	                    \end{center}

           \item[(d)] External associativity condition: The following
           diagram commutes.
\end{enumerate} \vspace{-.3cm}
 	            	          \begin{center}
      	          \resizebox{5in}{!}{ 
             $$
             \diagram ((U\otimes_j V)\otimes_j W)\otimes_i
            ((X\otimes_j Y)\otimes_j Z)
            \xto[rrr]^{\eta^{ij}_{(U\otimes_j V)W(X\otimes_j Y)Z}}
            \ar[d]^{\alpha^j \otimes_i \alpha^j} &&& \bigl((U\otimes_j
            V)\otimes_i (X\otimes_j Y)\bigr)\otimes_j(W\otimes_i Z)
            \dto^{\eta^{ij}_{UVXY}\otimes_j 1_{W\otimes_i Z}}\\
            (U\otimes_j (V\otimes_j W))\otimes_i (X\otimes_j
            (Y\otimes_j Z)) \dto^{\eta^{ij}_{U(V\otimes_j
            W)X(Y\otimes_j Z)}} &&&((U\otimes_i X)\otimes_j(V\otimes_i
            Y))\otimes_j(W\otimes_i Z) \ar[d]^{\alpha^j} \\
            (U\otimes_i X)\otimes_j\bigl((V\otimes_j W)\otimes_i
            (Y\otimes_j Z)\bigr) \xto[rrr]^{1_{U\otimes_i
            X}\otimes_j\eta^{ij}_{VWYZ}} &&& (U\otimes_i
            X)\otimes_j((V\otimes_i Y)\otimes_j(W\otimes_i Z))
            \enddiagram
            $$
          }
	           	                    \end{center}
\begin{enumerate}   
          \item[(e)] Finally it is required  for each triple $(i,j,k)$ satisfying
          $1\le i<j<k\le n$ that
          the giant hexagonal interchange diagram commutes.
         \end{enumerate}\vspace{-.8cm}
  		          \begin{center}
  	          \resizebox{5.3in}{!}{
          $$
          \xymatrix@C=-118pt{
          &((A\otimes_k A')\otimes_j (B\otimes_k B'))\otimes_i((C\otimes_k C')\otimes_j (D\otimes_k D'))
          \ar[ddl]|{\eta^{jk}_{AA'BB'}\otimes_i \eta^{jk}_{CC'DD'}}
          \ar[ddr]|{\eta^{ij}_{(A\otimes_k A')(B\otimes_k B')(C\otimes_k C')(D\otimes_k D')}}
          \\\\
          ((A\otimes_j B)\otimes_k (A'\otimes_j B'))\otimes_i((C\otimes_j D)\otimes_k (C'\otimes_j D'))
          \ar[dd]|{\eta^{ik}_{(A\otimes_j B)(A'\otimes_j B')(C\otimes_j D)(C'\otimes_j D')}}
          &&((A\otimes_k A')\otimes_i (C\otimes_k C'))\otimes_j((B\otimes_k B')\otimes_i (D\otimes_k D'))
          \ar[dd]|{\eta^{ik}_{AA'CC'}\otimes_j \eta^{ik}_{BB'DD'}}
          \\\\
          ((A\otimes_j B)\otimes_i (C\otimes_j D))\otimes_k((A'\otimes_j B')\otimes_i (C'\otimes_j D'))
          \ar[ddr]|{\eta^{ij}_{ABCD}\otimes_k \eta^{ij}_{A'B'C'D'}}
          &&((A\otimes_i C)\otimes_k (A'\otimes_i C'))\otimes_j((B\otimes_i D)\otimes_k (B'\otimes_i D'))
          \ar[ddl]|{\eta^{jk}_{(A\otimes_i C)(A'\otimes_i C')(B\otimes_i D)(B'\otimes_i D')}}
          \\\\
          &((A\otimes_i C)\otimes_j (B\otimes_i D))\otimes_k((A'\otimes_i C')\otimes_j (B'\otimes_i D'))
          }
          $$
          }
  		                    \end{center}

        \end{definition}

  \begin{definition} An {\it $n$--fold monoidal functor}
  $(F,\lambda^1,\dots,\lambda^n):{\cal C}\to{\cal D}$ between $n$--fold monoidal categories
  consists of a functor $F$ such that $F(I)=I$ together with natural
  transformations
  $$\lambda^i_{AB}:F(A)\otimes_i F(B)\to F(A\otimes_i B)\quad i=1,2,\dots, n$$
  satisfying the same associativity and unit conditions as monoidal functors.
  In addition the following hexagonal interchange diagram commutes.  
  \noindent
  	            	          \begin{center}
      	          \resizebox{5in}{!}{
    $$
  \diagram
  (F(A)\otimes_j F(B))\otimes_i(F(C)\otimes_j F(D))
  \xto[rrr]^{ {\eta^{ij}_{F(A)F(B)F(C)F(D)}}}
  \dto^{ {\lambda^j_{AB}\otimes_i\lambda^j_{CD}}}
  &&&(F(A)\otimes_i F(C))\otimes_j(F(B)\otimes_i F(D))
  \dto^{ \lambda^i_{AC}\otimes_j\lambda^i_{BD}}\\
  F(A\otimes_j B)\otimes_i F(C\otimes_j D)
  \dto^{ {\lambda^i_{(A\otimes_j B)(C\otimes_j D)}}}
  &&&F(A\otimes_i C)\otimes_j F(B\otimes_i D)
  \dto^{ {\lambda^j_{(A\otimes_i C)(B\otimes_i D)}}}\\
  F((A\otimes_j B)\otimes_i(C\otimes_j D))
  \xto[rrr]^{ F(\eta^{ij}_{ABCD})}
  &&& F((A\otimes_i C)\otimes_j(B\otimes_i D))
  \enddiagram
  $$ 
  }
         	                    \end{center}
 
  \end{definition}
  Composition of $n$-fold monoidal functors is defined 
  as for monoidal functors.

   The authors of \cite{Balt} point out that it is necessary to check that an $(n+1)$--fold monoidal category is
  the same thing as a tensor object in {\textbf{n}$\mathbf{ -MonCat}$, the category of
  $n$--fold monoidal categories and functors. Also as noticed in \cite{Balt}, the hexagonal interchange diagrams
  for the $(n+1)$--st monoidal operation regarded as an $n$--fold monoidal functor are what
  give rise to the giant hexagonal diagrams involving $\otimes_i$, $\otimes_j$ and $\otimes_{n+1}$.

  The authors of \cite{Balt}
  note that a symmetric monoidal category is $n$-fold monoidal for all $n$.  Just let
  $$\otimes_1=\otimes_2=\dots=\otimes_n=\otimes$$
  and define (associators added by myself)
  $$\eta^{ij}_{ABCD}=\alpha^{-1}\circ (1_A\otimes \alpha)\circ (1_A\otimes (c_{BC}\otimes 1_D))\circ (1_A\otimes \alpha^{-1})\circ \alpha$$
  for all $i<j$.

\section{Categories enriched over a $k$--fold monoidal category}

  \begin{theorem} \label{main:simple} For ${\cal V}$ a $k$--fold monoidal category ${\cal
  V}$--Cat is a $(k-1)$--fold
  monoidal 2-category. 
   \end{theorem} 
  
{\bf Example 1}

  We begin by describing the $k=2$ case. ${\cal V}$ is 2--fold
  monoidal with products $\otimes_{1}, \otimes_{2}.$ ${\cal
  V}$--categories (which are the objects of ${\cal V}$--Cat) are
  defined as being enriched over $({\cal
  V}$,$\otimes_{1},\alpha^1,I)$.  Here $\otimes_{1}$ plays the role of
  the product given by $\otimes$ in the axioms of section 1. We need
  to show that ${\cal V}$--Cat has a product.
  
  The unit object in ${\cal V}$--Cat is the enriched category ${\cal I}$ where $\left|{\cal I}\right| = \{0\}$ and 
  ${\cal I}(0,0) = I$. Of course $M_{000} = 1 = j_0.$
  The
  objects of the tensor ${\cal A} \otimes^{(1)}_{1}{\cal B}$ of two ${\cal
  V}$-categories
  ${\cal A}$ and ${\cal B}$ are simply pairs of objects, that is, elements
  of $\left|{\cal A}\right|\times\left|{\cal B}\right|$.   The hom--objects in ${\cal V}$ are given by
  $({\cal A} \otimes^{(1)}_{1}{\cal B})((A,B),(A',B')) = {\cal
  A}(A,A')\otimes_{2}{\cal B}(B,B')$. The composition morphisms that make
  ${\cal A} \otimes^{(1)}_{1}{
  \cal B}$ into a ${\cal V}$--category are
  immediately apparent as generalizations of the braided case. Recall that we are describing ${\cal
  A} \otimes^{(1)}_{1}{\cal B}$ as a 
  category enriched over ${\cal V}$ with product $\otimes_{1}$. Thus 
  \begin{small}
  $$M_{(A,B)(A',B')(A'',B'')} : ({\cal A} \otimes^{(1)}_{1}{\cal
  B})((A',B'),(A'',B''))\otimes_{1}({\cal A} \otimes^{(1)}_{1}{\cal
  B})((A,B),(A',B'))$$ $$\to ({\cal A} \otimes^{(1)}_{1}{\cal B})((A,B),(A'',B''))$$
  \end{small}
  is given by:
  $$
  \xymatrix{
  ({\cal A} \otimes^{(1)}_{1}{\cal B})((A',B'),(A'',B''))\otimes_{1}({\cal
  A} \otimes^{(1)}_{1}{\cal B})((A,B),(A',B'))
  \ar@{=}[d]\\
  ({\cal A}(A',A'')\otimes_{2}{\cal B}(B',B''))\otimes_{1}({\cal
  A}(A,A')\otimes_{2}{\cal B}(B,B'))
  \ar[d]_{\eta^{1,2}}\\
  ({\cal A}(A',A'')\otimes_{1}{\cal A}(A,A'))\otimes_{2}({\cal B}(B',B'')\otimes_{1}{\cal B}(B,B'))
  \ar[d]_{M_{AA'A''}\otimes_{2}M_{BB'B''}}\\
  ({\cal A}(A,A'')\otimes_{2}{\cal B}(B,B''))
  \ar@{=}[d]\\
  ({\cal A} \otimes^{(1)}_{1}{\cal B})((A,B),(A'',B''))
  }
  $$

{\bf Example 2}

  Next we describe the $k=3$ case. ${\cal V}$ is 3--fold monoidal with
  products $\otimes_{1}, \otimes_{2}$ and $\otimes_{3}$. ${\cal
  V}$--categories are defined as being enriched over $({\cal
  V}$,$\otimes_{1},\alpha^1,I).$ Now ${\cal V}$--Cat has two
  products. The objects of both possible tensors ${\cal A}
  \otimes^{(1)}_{1}{\cal B}$ and ${\cal A} \otimes^{(1)}_{2}{\cal B}$
  of two ${\cal V}$-categories ${\cal A}$ and ${\cal B}$ are elements
  in $\left|{\cal A}\right|\times\left|{\cal B}\right|$.  The
  hom--objects in ${\cal V}$ are given by
  $$({\cal A} \otimes^{(1)}_{1}{\cal B})((A,B),(A',B')) = {\cal A}(A,A')\otimes_{2}{\cal B}(B,B')$$
  just as in the previous case, and by
  $$({\cal A} \otimes^{(1)}_{2}{\cal B})((A,B),(A',B')) = {\cal
  A}(A,A')\otimes_{3}{\cal B}(B,B').$$ 
  The composition that makes $({\cal A} \otimes^{(1)}_{2}{\cal B})$ into
  a ${\cal V}$--category is analogous to that for $({\cal A}
  \otimes^{(1)}_{1}{\cal B})$ but uses $\eta^{1,3}$ as its middle
  exchange morphism.
  
  Now we need an interchange 2--natural transformation $\eta^{(1)1,2}$ for ${\cal V}$--Cat. 
  The family of morphisms $\eta^{(1)1,2}_{{\cal A}{\cal B}{\cal C}{\cal D}}$
  that make up a 2--natural transformation between the 2--functors  
  $\times^4{\cal V}$--Cat $:\to {\cal V}$--Cat in question 
  is a family of enriched functors. 
  Their action on objects
  is to send
  $$((A,B),(C,D)) \in \left|({\cal A} \otimes^{(1)}_{2}{\cal B}) \otimes^{(1)}_{1}({\cal C} \otimes^{(1)}_{2}{\cal D})\right|
   $$ $$\text{ to } ((A,C),(B,D)) \in \left|({\cal A} \otimes^{(1)}_{1}{\cal C}) \otimes^{(1)}_{2}({\cal B} \otimes^{(1)}_{1}{\cal D})\right|.$$
  The correct construction of the family of hom--object morphisms in ${\cal V}$--Cat for each of these functors 
  is also clear.
  Noting that 
  $$[({\cal A} \otimes^{(1)}_{2} {\cal B}) \otimes^{(1)}_{1}({\cal C} \otimes^{(1)}_{2}{\cal D})](((A,B),(C,D)),((A',B'),(C',D')))$$
  $$=({\cal A} \otimes^{(1)}_{2}{\cal B})((A,B),(A',B'))\otimes_{2}({\cal C} \otimes^{(1)}_{2}{\cal D})((C,D),(C',D')) $$
  $$=({\cal A}(A,A')\otimes_{3}{\cal B}(B,B'))\otimes_{2}({\cal C}(C,C')\otimes_{3}{\cal D}(D,D'))$$
  and similarly 
  $$[({\cal A} \otimes^{(1)}_{1}{\cal C}) \otimes^{(1)}_{2}({\cal B} \otimes^{(1)}_{1}{\cal D})](((A,C),(B,D)),((A',C'),(B',D')))$$
  $$=({\cal A}(A,A')\otimes_{2}{\cal C}(C,C'))\otimes_{3}({\cal B}(B,B')\otimes_{2}{\cal D}(D,D'))$$
  we make the obvious identification, where by obvious I mean based upon the corresponding structure in 
  ${\cal V}.$ For a detailed discussion of this construction for the case of braided ${\cal V}$ see \cite{forcey1}.
 Here ``based upon'' is more freely 
  interpreted as allowing a shift in index.
  Thus we write:
  $$\eta^{(1)1,2}_{{{\cal A}{\cal B}{\cal C}{\cal
  D}}_{(ABCD)(A'B'C'D')}} = \eta^{2,3}_{{\cal A}(A,A'){\cal
  B}(B,B'){\cal C}(C,C'){\cal D}(D,D')}$$ Much needs to be
  verified. Existence and coherence of required natural
  transformations, satisfaction of enriched axioms and of $k$--fold
  monoidal axioms all must be checked. These will be dealt with next.
  
  \begin{proof}[Proof of Theorem \ref{main:simple}]As in the examples,
  ${\cal V}$--Cat is made up of categories enriched over $({\cal
  V}$,$\otimes_{1},\alpha^1,I).$ Here we define products $
  \otimes^{(1)}_{1} \ldots \otimes^{(1)}_{k-1}$ in ${\cal V}$--Cat for
  ${\cal V}$ $k$--fold monoidal.  We check that our products do make
  ${\cal A} \otimes^{(1)}_{2}{\cal B}$ into a ${\cal V}$--category.
  Then we check that ${\cal V}$--Cat has the required coherent
  2--natural transformations of associativity and units. We then
  define interchange 2--natural transformations $\eta^{(1)i,j}$ and
  check that the interchange transformations are 2--natural and obey
  all the axioms required of them. It is informative to observe how
  these axioms are satisfied based upon the axioms that ${\cal V}$
  itself satisfies. It is here that we should look carefully for the
  algebraic reflection of the topological functor $\Omega.$
  
  Again, the unit object in ${\cal V}$--Cat is the enriched category ${\cal I}$ where $\left|{\cal I}\right| = \{0\}$ and 
  ${\cal I}(0,0) = I$.
  For ${\cal V}$ $k$--fold monoidal we define the $i$th product of ${\cal
  V}$--categories ${\cal A} \otimes^{(1)}_{i}{\cal B}$
   to have objects $\in \left|{\cal A}\right|\times \left|{\cal B}\right|$
   and to have hom--objects in ${\cal V}$ given by
   $$({\cal A} \otimes^{(1)}_{i} {\cal B})((A,B),(A',B')) = {\cal
   A}(A,A')\otimes_{i+1} {\cal B}(B,B').$$ 
   Immediately we see that ${\cal V}$--Cat is $(k-1)$--fold monoidal by
   definition.  The composition morphisms are \begin{small}
   $$M_{(A,B)(A',B')(A'',B'')} : ({\cal A} \otimes^{(1)}_{i}{\cal
   B})((A',B'),(A'',B''))\otimes_{1}({\cal A} \otimes^{(1)}_{i}{\cal
   B})((A,B),(A',B'))$$ $$\to ({\cal A} \otimes^{(1)}_{i}{\cal B})((A,B),(A'',B''))$$
   \end{small}
   given by
   $$
   \xymatrix{
   ({\cal A} \otimes^{(1)}_{i}{\cal B})((A',B'),(A'',B''))\otimes_{1}({\cal
   A} \otimes^{(1)}_{i}{\cal B})((A,B),(A',B'))
   \ar@{=}[d]\\
   ({\cal A}(A',A'')\otimes_{i+1}{\cal B}(B',B''))\otimes_{1}({\cal
   A}(A,A')\otimes_{i+1}{\cal B}(B,B'))
   \ar[d]_{\eta^{1,i+1}}\\
   ({\cal A}(A',A'')\otimes_{1}{\cal A}(A,A'))\otimes_{i+1}({\cal
   B}(B',B")\otimes_{1}{\cal B}(B,B'))
   \ar[d]_{M_{AA'A''}\otimes_{2}M_{BB'B''}}\\
   ({\cal A}(A,A'')\otimes_{i+1}{\cal B}(B,B''))
   \ar@{=}[d]\\
   ({\cal A} \otimes^{(1)}_{i}{\cal B})((A,B),(A'',B'')).
  }
  $$
   The identity element is given by $j_{(A,B)} = \xymatrix{I = I \otimes_{i+1} I
   							\ar[d]^{j_A \otimes_{i+1} j_B}
   							\\{\cal A}(A,A)\otimes_{i+1} {\cal B}(B,B)
   							\ar@{=}[d]
   							\\({\cal A} \otimes^{(1)}_{i}{\cal B})((A,B),(A,B)).
   							}$
   
   The product $\otimes^{(1)}_i$ of enriched functors is defined in the obvious way.
   
  Here we first check that ${\cal A}\otimes^{(1)}_{i}{\cal B}$ is indeed properly
   enriched over ${\cal V}.$ Our definition of $M$ must obey the axioms for associativity and respect of the unit.
   For associativity the following diagram must commute, where the initial bullet represents
   $$
  [({\cal A}\otimes^{(1)}_{i}{\cal B})((A'',B''),(A''',B'''))\otimes_{1} ({\cal A}\otimes^{(1)}_{i}{\cal B})((A',B'),(A'',B''))]
  $$ $$\otimes_{1}({\cal A}\otimes^{(1)}_{i}{\cal B})((A,B),(A',B')).
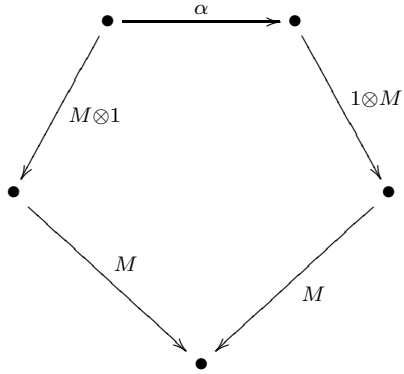
  $$
    $$
    \xymatrix{
    &\bullet
    \ar[rr]^{ \alpha}
    \ar[ddl]^{ M \otimes 1}
    &&\bullet
    \ar[ddr]^{ 1 \otimes M}&\\\\
    \bullet
    \ar[ddrr]^{ M}
    &&&&\bullet
    \ar[ddll]^{ M}
    \\\\&&\bullet
    }$$
  In the expanded diagram given in Figure \ref{diag1} let $X={\cal
  A}(A,A')$, $X'={\cal A}(A',A'')$, $X''={\cal A}(A'',A''')$, $Y={\cal
  B}(B,B')$, $Y'={\cal B}(B',B'')$ and $Y''={\cal B}(B'',B''').$ The
  exterior of the diagram is required to commute.

\begin{figure}[p] 
\begin{center}
    	 \resizebox{4.8in}{!}{
      \begin{sideways}
      \begin{small}
      $$
      \xymatrix@C=-95pt@R=30pt{
      && *+<10pt>{\text{ }}
      \ar[ddl]^{\eta^{1,i+1}\otimes_{1} 1}
      &
      *\txt{$[(X''\otimes_{i+1} Y'')\otimes_{1} (X'\otimes_{i+1} Y')]\text{ }\text{ }\text{ }\text{ }\text{ }\text{ }\text{ }\text{ }\text{ }\text{ }\text{ }\text{ }\text{ }\text{ }\text{ }\text{ }\text{ }\text{ }\text{ }\text{ }\text{ }\text{ }\text{ }\text{ }\text{ }\text{ }\text{ }\text{ }\text{ }\text{ }\text{ }\text{ }\text{ }\text{ }\text{ }\text{ }\text{ }\text{ }\text{ }\text{ }\text{ }\text{ }\text{ }\text{ }\text{ }\text{ }$
      \\$\otimes_{1} (X\otimes_{i+1} Y){-}^{\alpha^{1}}\text{ \normalsize{$\to$}}(X''\otimes_{i+1} Y'')$
      \\$\text{ }\text{ }\text{ }\text{ }\text{ }\text{ }\text{ }\text{ }\text{ }\text{ }\text{ }\text{ }\text{ }\text{ }\text{ }\text{ }\text{ }\text{ }\text{ }\text{ }\text{ }\text{ }\text{ }\text{ }\text{ }\text{ }\text{ }\text{ }\text{ }\text{ }\text{ }\text{ }\text{ }\text{ }\text{ }\text{ }\text{ }\text{ }\text{ }\text{ }\text{ }\text{ }\text{ }\text{ }\text{ }\text{ }\text{ }\otimes_{1} [(X'\otimes_{i+1} Y')\otimes_{1} (X\otimes_{i+1} Y)]$}
      &*+<45pt>{\text{ }}
      \ar[ddr]^{1 \otimes_{1} \eta^{1,i+1}}
      \\\\
      &[(X''\otimes_{1} X')\otimes_{i+1} (Y''\otimes_{1} Y')]\otimes_{1} (X\otimes_{i+1} Y)
      \ar[dddl]_{(M\otimes_{i+1} M)\otimes_{1} (1\otimes_{i+1} 1)\text{  }}
      \ar[ddr]^{\eta^{1,i+1}}
      &&&&(X''\otimes_{i+1} Y'')\otimes_{1} [(X'\otimes_{1} X)\otimes_{i+1} (Y'\otimes_{1} Y)]
      \ar[ddl]^{\eta^{1,i+1}}
      \ar[dddr]^{(1\otimes_{i+1} 1)\otimes_{1} (M\otimes_{i+1} M)}
      \\\\
      && *+<50pt>{\text{ }}\text{ }\text{ }\text{ }\text{ }\text{ }\text{ } \text{ }\text{ }\text{ }\text{ }\text{ }\text{ }
      \ar[ddl]^{\text{  }(M\otimes_{1} 1)\otimes_{i+1} (M\otimes_{1} 1)}
      &
      *\txt{$[(X''\otimes_{1} X')\otimes_{1} X]\text{ }\text{ }\text{ }\text{ }\text{ }\text{ }\text{ }\text{ }\text{ }\text{ }\text{ }\text{ }\text{ }\text{ }\text{ }\text{ }\text{ }\text{ }\text{ }\text{ }\text{ }\text{ }\text{ }\text{ }\text{ }\text{ }\text{ }\text{ }\text{ }\text{ }\text{ }\text{ }\text{ }\text{ }\text{ }\text{ }\text{ }\text{ }\text{ }\text{ }\text{ }\text{ }\text{ }\text{ }$
      \\$\otimes_{i+1} [(Y''\otimes_{1} Y')\otimes_{1} Y]{-}^{\alpha^{1}\otimes_{i+1} \alpha^{1}}\text{ \normalsize{$\to$}}[X''\otimes_{1} (X'\otimes_{1} X)]$
      \\$\text{ }\text{ }\text{ }\text{ }\text{ }\text{ }\text{ }\text{ }\text{ }\text{ }\text{ }\text{ }\text{ }\text{ }\text{ }\text{ }\text{ }\text{ }\text{ }\text{ }\text{ }\text{ }\text{ }\text{ }\text{ }\text{ }\text{ }\text{ }\text{ }\text{ }\text{ }\text{ }\text{ }\text{ }\text{ }\text{ }\text{ }\text{ }\text{ }\text{ }\text{ }\text{ }\text{ }\text{ }\text{ }\text{ }\text{ }\text{ }\text{ }\text{ }\text{ }\text{ }\otimes_{i+1} [Y''\otimes_{1} (Y'\otimes_{1} Y)]$}
      &*+<50pt>{\text{ }}\text{ }\text{ }\text{ }\text{ }\text{ }\text{ }\text{ }\text{ }\text{ }
      \ar[ddr]_{(1\otimes_{1} M)\otimes_{i+1} (1\otimes_{1} M)\text{  }}
      &&
      \\
      ({\cal A}(A',A''')\otimes_{i+1} {\cal B}(B',B'''))\otimes_{1} (X\otimes_{i+1} Y)
      \ar[dr]^{\eta^{1,i+1}}
      &&&&&&
      (X''\otimes_{i+1} Y'')\otimes_{1} ({\cal A}(A,A'')\otimes_{i+1} {\cal B}(B,B''))
      \ar[dl]^{\eta^{1,i+1}}
      \\
      &({\cal A}(A',A''')\otimes_{1} X)\otimes_{i+1} ({\cal B}(B',B''')\otimes_{1} Y)
      \ar[ddrr]|{M\otimes_{i+1} M}
      &&&&(X''\otimes_{1} {\cal A}(A,A''))\otimes_{i+1} (Y''\otimes_{1} {\cal B}(B,B''))
      \ar[ddll]|{M\otimes_{i+1} M}
      \\\\
      &&&{\cal A}(A,A''')\otimes_{i+1} {\cal B}(B,B''')
      } 
      $$
      \end{small}
      \end{sideways}
       }
  \end{center}\vspace{-1cm}
  \nocolon\caption{}\label{diag1}
 \end{figure}

  The lower pentagon in Figure \ref{diag1} commutes since it is two
  copies of the associativity axiom--one for ${\cal A}$ and one for
  ${\cal B}.$ The two diamonds commute by the naturality of $\eta.$
  The upper hexagon commutes by the internal associativity of $\eta.$

 For the unit axioms we have the following compact diagram.
              	\vspace{-8mm}	          \begin{center}
	          \resizebox{5.5in}{!}
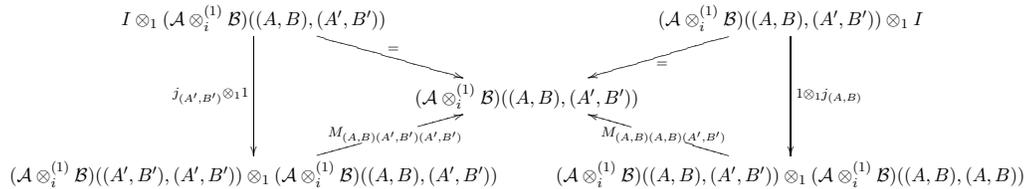
{
  $$
      \xymatrix@C=-30pt{
      I\otimes_1 ({\cal A} \otimes^{(1)}_{i}{\cal B})((A,B),(A',B'))
      \ar[rrd]^{=}
      \ar[dd]_{j_{(A',B')}\otimes_1 1}
      &&&&({\cal A} \otimes^{(1)}_{i}{\cal B})((A,B),(A',B'))\otimes_1 I 
      \ar[dd]^{1\otimes_1 j_{(A,B)}}
      \ar[lld]^{=}\\
      &&({\cal A} \otimes^{(1)}_{i}{\cal B})((A,B),(A',B'))\\
      ({\cal A} \otimes^{(1)}_{i}{\cal B})((A',B'),(A',B'))\otimes_1 ({\cal A} \otimes^{(1)}_{i}{\cal B})((A,B),(A',B'))
      \ar[rru]|{\text{ }\text{ }M_{(A,B)(A',B')(A',B')}}
      &&&&({\cal A} \otimes^{(1)}_{i}{\cal B})((A,B),(A',B'))\otimes_1 ({\cal A} \otimes^{(1)}_{i}{\cal B})((A,B),(A,B))
      \ar[llu]|{\text{ }\text{ }M_{(A,B)(A,B)(A',B')}}
      }
   $$ 
  } \end{center} 
  I expand the left triangle, abbreviating $X={\cal
  A}(A,A')$, $Y={\cal A}(A',A')$, $Z={\cal B}(B,B')$ and $W={\cal
  B}(B',B').$ The exterior of the following must commute:
  \begin{center} \resizebox{4.5in}{!}
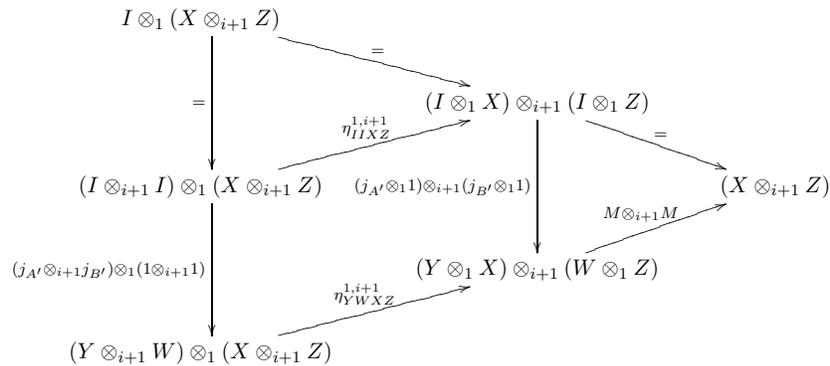
{
      $$
  \xymatrix{
  I\otimes_1 (X\otimes_{i+1} Z)\text{ }\text{ }\text{ }
  \ar[dd]_{=}
  \ar[dr]^{=}
  \\
  &(I\otimes_1 X)\otimes_{i+1} (I\otimes_1 Z)
  \ar[dr]^{=}
  \ar[dd]_{(j_{A'} \otimes_1 1)\otimes_{i+1} (j_{B'}\otimes_1 1)}
  \\
  (I\otimes_{i+1} I)\otimes_1 (X\otimes_{i+1} Z)\text{ }\text{ }\text{ }
  \ar[dd]_{(j_{A'} \otimes_{i+1} j_{B'})\otimes_1 (1\otimes_{i+1} 1)}
  \ar[ur]^{\eta^{1,i+1}_{IIXZ}}
  &&(X\otimes_{i+1} Z)
  \\
  &(Y\otimes_1 X)\otimes_{i+1} (W\otimes_1 Z)
  \ar[ur]^{M\otimes_{i+1} M}
  \\
  (Y\otimes_{i+1} W)\otimes_1 (X\otimes_{i+1} Z)\text{ }\text{ }\text{ }
  \ar[ur]^{\eta^{1,i+1}_{YWXZ}}
  }
  $$
  } \end{center} 
      The parallelogram commutes by naturality of $\eta$, the rightmost
      triangle by the unit axioms of the individual ${\cal
      V}$--categories, and the top triangle by the internal unit
      condition for $\eta.$ The right-hand triangle in the axiom is
      checked similarly.
  
  On a related note, we need to check that ${\cal I} \otimes^{(1)}_i
  {\cal A}$ = ${\cal A}$ The object sets and hom--objects of the two
  categories in question are clearly equivalent. What needs to be
  checked is that the composition morphisms are the same.  Note that
  the composition given by
  $$
   \xymatrix{
   ({\cal I} \otimes^{(1)}_{i}{\cal A})((0,A'),(0,A''))\otimes_{1} ({\cal I} \otimes^{(1)}_{i} {\cal A})((0,A),(0,A'))
   \ar@{=}[d]\\
   (I\otimes_{i+1} {\cal A}(A',A''))\otimes_{1} (I\otimes_{i+1} {\cal A}(A,A'))
   \ar[d]_{\eta^{1,i+1}_{I{\cal A}(A',A'')I{\cal A}(A,A')}}\\
   (I\otimes_{1}I)\otimes_{i+1}({\cal A}(A',A'')\otimes_{1}{\cal A}(A,A'))
   \ar[d]_{1\otimes_{i+t} M_{AA'A''}}\\
   (I\otimes_{i+1}{\cal A}(A,A''))
   \ar@{=}[d]\\
   ({\cal I} \otimes^{(1)}_{i} {\cal A})((0,A),(0,A''))
  }
  $$
  is equivalent to simply $M_{AA'A''}$ by the external unit condition for $\eta.$

   Associativity in ${\cal V}$--Cat
   must hold for each $\otimes^{(1)}_{i}$. The components of 2--natural
   isomorphism 
   $$\alpha^{(1)i}_{{\cal A}{\cal B}{\cal C}}: ({\cal A} \otimes^{(1)}_{i} {\cal B})
    \otimes^{(1)}_{i} {\cal C} \to {\cal A} \otimes^{(1)}_{i} ({\cal B} \otimes^{(1)}_{i}
   {\cal C})$$
   are ${\cal V}$--functors
   that send ((A,B),C) to (A,(B,C)) and whose hom-components 
   \begin{small}
   $$\alpha^{(1)i}_{{{\cal A}{\cal B}{\cal C}}_{((A,B),C)((A',B'),C')}}: [({\cal A} \otimes^{(1)}_{i} {\cal B}) \otimes^{(1)}_{i} {\cal C}](((A,B),C),((A',B'),C'))
   $$ $$\to [{\cal A} \otimes^{(1)}_{i} ({\cal B} \otimes^{(1)}_{i} {\cal C})]((A,(B,C)),(A',(B',C')))$$
   \end{small}
   are given by: 
   $$\alpha^{(1)i}_{{{\cal A}{\cal B}{\cal C}}_{((A,B),C)((A',B'),C')}}
   = \alpha^{i+1}_{{\cal A}(A,A'){\cal B}(B,B'){\cal C}(C,C')}$$
   This guarantees that the 2--natural isomorphism $\alpha^{(1)i}$ is coherent. The commutativity of the pentagon for
   the objects is trivial, and the commutativity of the pentagon for the hom--object morphisms follows directly from the 
   commutativity of the pentagon for $\alpha^{i+1}.$
  
   In order to be a functor the associator components must
   satisfy the commutativity of the
   diagrams in Definition \ref{enriched:funct}.
   $$
     \xymatrix{
     &\bullet
     \ar[rr]^{M}
     \ar[d]^{\alpha^{(1)i} \otimes \alpha^{(1)i}}
     &&\bullet
     \ar[d]^{\alpha^{(1)i}}&\\
     &\bullet
     \ar[rr]^{M}
     &&\bullet
     }\leqno{(1)}
   $$
   $$
     \xymatrix{
     &&\bullet
     \ar[dd]^{ \alpha^{(1)i}}\\
     I
     \ar[rru]^{j_{((A,B),C)}}
     \ar[rrd]_{ j_{(A,(B,C))}}\\
     &&\bullet
     }\leqno{(2)}
   $$
   Expanding the first using
   the definitions just given we have that the initial position in the diagram is 
\vspace{-2mm}

\resizebox{4.5in}{!}{\vbox{
   $$
   [({\cal A}\otimes^{(1)}_{i}{\cal B})\otimes^{(1)}_{i}{\cal C}](((A',B'),C'),((A'',B''),C''))\otimes_{1}[({\cal A}\otimes^{(1)}_{i}{\cal B})\otimes^{(1)}_{i}{\cal C}](((A,B),C),((A',B'),C'))
  $$ $$
   = [({\cal A}(A',A'')\otimes_{i+1}{\cal B}(B',B''))\otimes_{i+1}{\cal C}(C',C'')]\otimes_{1}[({\cal A}(A,A')\otimes_{i+1}{\cal B}(B,B'))\otimes_{i+1}{\cal C}(C,C')]
   $$}
}
   
\vspace{-2mm}
We let $X={\cal A}(A',A'')$, $Y={\cal B}(B',B'')$, $Z={\cal C}(C',C'')$, $X'={\cal A}(A,A')$, $Y'={\cal B}(B,B')$ and $Z'={\cal C}(C,C').$
   Then expanding the diagram, with an added interior arrow, we have:
               		          \begin{center}
	          \resizebox{5in}{!}{
   $$
   \xymatrix@C=-100pt{
   &[(X\otimes_{i+1} Y)\otimes_{i+1} Z]\otimes_{1} [(X'\otimes_{i+1} Y')\otimes_{i+1} Z']
   \ar[ddr]|{\text{  }\eta^{1,i+1}_{(X\otimes_{i+1} Y)Z(X'\otimes_{i+1} Y')Z'}}
   \ar[ddl]|{\alpha^{i+1}\otimes_{1} \alpha^{i+1}}\\\\
   [X\otimes_{i+1} (Y\otimes_{i+1} Z)]\otimes_{1} [X'\otimes_{i+1} (Y'\otimes_{i+1} Z')]
   \ar[dd]|{\eta^{1,i+1}_{X(Y\otimes_{i+1} Z)X'(Y'\otimes_{i+1} Z')}}
   &&[(X\otimes_{i+1} Y)\otimes_{1} (X'\otimes_{i+1} Y')]\otimes_{i+1} (Z\otimes_{1} Z')
   \ar[dd]|{\eta^{1,i+1}_{XYX'Y'}\otimes_{i+1} 1_{Z\otimes_{1} Z'}}\\\\
   (X\otimes_{1} X') \otimes_{i+1} [(Y\otimes_{i+1} Z)\otimes_{1} (Y'\otimes_{i+1} Z')]
   \ar[dd]|{1_{X\otimes_{1} X'} \otimes_{i+1} \eta^{1,i+1}_{YZY'Z'}}
   &&[(X\otimes_{1} X')\otimes_{i+1} (Y\otimes_{1} Y')]\otimes_{i+1} (Z\otimes_{1} Z')
   \ar[dd]|{(M\otimes_{i+1} M)\otimes_{i+1} M}
   \ar[ddll]|{\alpha^{i+1}}
   \\\\
   (X\otimes_{1} X')\otimes_{i+1} [(Y\otimes_{1} Y')\otimes_{i+1} (Z\otimes_{1} Z')]
   \ar[ddr]|{M\otimes_{i+1} (M\otimes_{i+1} M)}
   &&({\cal A}(A,A'')\otimes_{i+1} {\cal B}(B,B''))\otimes_{i+1}{\cal C}(C,C'')
   \ar[ddl]|{\alpha^{i+1}}\\\\
   &{\cal A}(A,A'')\otimes_{i+1} ({\cal B}(B,B'')\otimes_{i+1}{\cal C}(C,C''))
   }
   $$
   }
   \end{center}

   The lower quadrilateral    commutes by
   naturality of $\alpha$, and the upper hexagon commutes by the
   external associativity of $\eta.$
   
    The uppermost position in the expanded version of diagram number (2) is
   $$
   [({\cal A}\otimes^{(1)}_{i}{\cal B})\otimes^{(1)}_{i}{\cal C}](((A,B),C),((A,B),C))
   $$
   $$
   = [({\cal A}(A,A)\otimes_{i+1}{\cal B}(B,B))\otimes_{i+1}{\cal C}(C,C)].
   $$
   The expanded diagram is easily seen to commute by the naturality of $\alpha.$
   
   The 2--naturality of $\alpha^{(1)}$ is essentially just the naturality of its 
   components, but I think it ought to be expounded upon. Since the components of $\alpha^{(1)}$ are
   ${\cal V}$--functors the whisker diagrams for the definition of 2--naturality are defined by the 
   whiskering in ${\cal V}$--Cat. Given an arbitrary 2--cell in $\times^3{\cal V}$--Cat, i.e. 
   $(\beta,\gamma,\rho):(Q,R,S)\to (Q',R',S'):({\cal A},{\cal B},{\cal C})\to({\cal A'},{\cal B'},{\cal C'})$
   the diagrams whose composition must be equal are:
  \vspace{-4mm}
                  \begin{center}
      	 \resizebox{5in}{!}{
         $$
    \xymatrix@R-=3pt{
    &\ar@{=>}[dd]^{(\beta\otimes^{(1)}_i \gamma)\otimes^{(1)}_i \rho}\\
     ({\cal A}\otimes^{(1)}_i {\cal B})\otimes^{(1)}_i {\cal C}
    \ar@/^2pc/[rrr]^{(Q\otimes^{(1)}_i R)\otimes^{(1)}_i S}
    \ar@/_2pc/[rrr]_{(Q'\otimes^{(1)}_i R')\otimes^{(1)}_i S'}
    &&&({\cal A}'\otimes^{(1)}_i {\cal B}')\otimes^{(1)}_i {\cal C}'
    \ar[rr]^{\alpha^{(1)i}_{{\cal A}'{\cal B}'{\cal C}'}}
    &&{\cal A}'\otimes^{(1)}_i ({\cal B}'\otimes^{(1)}_i {\cal C}')
    \\
   &\\
   }
   $$
   }
       \end{center}\vspace{-4mm}
       
                       \begin{center}
           	 \resizebox{5in}{!}{
             $$
  = \xymatrix@R-=3pt{
   &&&\ar@{=>}[dd]^{\beta\otimes^{(1)}_i (\gamma\otimes^{(1)}_i \rho)}\\
   ({\cal A}\otimes^{(1)}_i {\cal B})\otimes^{(1)}_i {\cal C}
   \ar[rr]^{\alpha^{(1)i}_{{\cal A}{\cal B}{\cal C}}}
   &&{\cal A}\otimes^{(1)}_i ({\cal B}\otimes^{(1)}_i {\cal C})
   \ar@/^2pc/[rrr]^{Q\otimes^{(1)}_i (R\otimes^{(1)}_i S)}
   \ar@/_2pc/[rrr]_{Q'\otimes^{(1)}_i (R'\otimes^{(1)}_i S')}
   &&&{\cal A}'\otimes^{(1)}_i ({\cal B}'\otimes^{(1)}_i {\cal C}')
   \\
   &&&&\\
   }
  $$ 
  }
    \end{center}
  This is quickly seen to hold when we translate using the definitions of whiskering in ${\cal V}$--Cat, as follows.
  The $ABCD$ components of the new 2--cells are given by the exterior legs of the following diagram. They are equal 
  by naturality of $\alpha^{i+1}$ and Mac Lane's coherence theorem.
                		\vspace{-4mm}         \begin{center}
	          \resizebox{5in}{!}{
  $$
  \xymatrix{
  &I
  \ar[dr]^{=}
  \ar[dl]^{=}
  \\
  (I\otimes_{i+1} I)\otimes_{i+1} I
  \ar[rr]^{\alpha^{i+1}}
  \ar[dd]^{(\beta_A\otimes_{i+1} \gamma_B\otimes_{i+1}) \rho_C}
  &&I\otimes_{i+1} (I\otimes_{i+1} I)
  \ar[dd]^{\beta_A\otimes_{i+1} (\gamma_B\otimes_{i+1} \rho_C)}
  \\\\
  ({\cal A}'(QA,Q'A)\otimes_{i+1} {\cal B}'(RB,R'B))\otimes_{i+1} {\cal C}'(SC,S'C)
  \ar[rr]^{\alpha^{i+1}}
  &\text{ }&{\cal A}'(QA,Q'A)\otimes_{i+1} ({\cal B}'(RB,R'B)\otimes_{i+1} {\cal C}'(SC,S'C))
  }
  $$
  }
  \end{center}
  
  Now we turn to consider the existence and behavior of interchange 2--natural transformations $\eta^{(1)ij}$
  for $j\ge i+1$.
  As in the example, we define the component morphisms $\eta^{(1)i,j}_{{\cal A}{\cal B}{\cal C}{\cal D}}$
  that make a 2--natural transformation between 2--functors. Each component must be an enriched functor.
  Their action on objects
  is to send $((A,B),(C,D)) \in \left|({\cal A}\otimes^{(1)}_{j} {\cal B})\otimes^{(1)}_{i} ({\cal C}\otimes^{(1)}_{j} {\cal D})\right|$
   to $((A,C),(B,D)) \in \left|({\cal A}\otimes^{(1)}_{i} {\cal C})\otimes^{(1)}_{j} ({\cal B}\otimes^{(1)}_{i} {\cal D})\right|$.
  The hom--object morphisms are given by:
   $$\eta^{(1)i,j}_{{{\cal A}{\cal B}{\cal C}{\cal D}}_{(ABCD)(A'B'C'D')}} =
   \eta^{i+1,j+1}_{{\cal A}(A,A'){\cal B}(B,B'){\cal C}(C,C'){\cal D}(D,D')}$$
  
  For this designation of $\eta^{(1)}$ to define a valid ${\cal V}$--functor, it must obey the axioms for compatibility 
  with composition and units. 
  We need commutativity of the following diagram, where the first bullet represents
  \begin{small}
  $$[({\cal A}\otimes^{(1)}_{j} {\cal B})\otimes^{(1)}_{i} ({\cal C}\otimes^{(1)}_{j} {\cal D})](((A',B'),(C',D')),((A'',B''),(C'',D'')))$$ $$\otimes_1 [({\cal A}\otimes^{(1)}_{j} {\cal B})\otimes^{(1)}_{i} ({\cal C}\otimes^{(1)}_{j} {\cal D})](((A,B),(C,D)),((A',B'),(C',D')))$$
  \end{small}
  and the last bullet represents
  \begin{small}
  $$[({\cal A}\otimes^{(1)}_{i} {\cal C})\otimes^{(1)}_{j} ({\cal B}\otimes^{(1)}_{i} {\cal D})](((A,C),(B,D)),((A'',C''),(B'',D''))).$$
  \end{small}
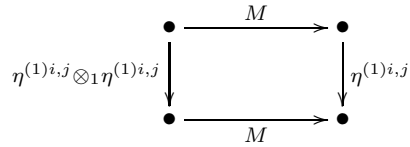
\vspace{-4mm}
  $$
  \xymatrix{
  \bullet
  \ar[rr]^{M}
  \ar[d]_{\eta^{(1)i,j} \otimes_1 \eta^{(1)i,j}}
  &&\bullet
  \ar[d]^{\eta^{(1)i,j}}
 \\
 \bullet
  \ar[rr]_{M}
  &&\bullet
  }
  $$
  If we let $X={\cal A}(A,A')$, $Y={\cal B}(B,B')$, $Z={\cal
  C}(C,C')$, $W={\cal D}(D,D')$, $X'={\cal A}(A',A'')$, $Y'={\cal
  B}(B',B'')$, $Z'={\cal C}(C',C'')$ and $W'={\cal D}(D',D'')$ then
  the expanded diagram is given in Figure \ref{diag2}.  The exterior
  must commute.

\begin{figure}[p]
\begin{center}
	          \resizebox{4.5in}{!}{
  $$
  \begin{sideways}
  \begin{small}
  \xymatrix@C=-140pt@R=30pt{
  \text{ }
  \\
  &[(X'\otimes_{j+1} Y')\otimes_{i+1} (Z'\otimes_{j+1} W')]\otimes_{1} [(X\otimes_{j+1} Y)\otimes_{i+1} (Z\otimes_{j+1} W)]
  \ar[ddr]|{\eta^{1,i+1}_{(X'\otimes_{j+1} Y')(Z'\otimes_{j+1} W')(X\otimes_{j+1} Y)(Z\otimes_{j+1} W)}}
  \ar[ddl]|{\eta^{i+1,j+1}_{X'Y'Z'W'}\otimes_1 \eta^{i+1,j+1}_{XYZW}}\\\\
  [(X'\otimes_{i+1} Z')\otimes_{j+1} (Y'\otimes_{i+1} W')]\otimes_{1} [(X\otimes_{i+1} Z)\otimes_{j+1} (Y\otimes_{i+1} W)]
  \ar[dd]|{\eta^{1,j+1}_{(X'\otimes_{i+1} Z')(Y'\otimes_{i+1} W')(X\otimes_{i+1} Z)(Y\otimes_{i+1} W)}}
  &&[(X'\otimes_{j+1} Y')\otimes_{1} (X\otimes_{j+1} Y)]\otimes_{i+1} [(Z'\otimes_{j+1} W')\otimes_{1} (Z\otimes_{j+1} W)]
  \ar[dd]|{\eta^{1,j+1}_{X'Y'XY}\otimes_{i+1} \eta^{1,j+1}_{Z'W'ZW}}\\\\
  [(X'\otimes_{i+1} Z')\otimes_{1} (X\otimes_{i+1} Z)]\otimes_{j+1} [(Y'\otimes_{i+1} W')\otimes_{1} (Y\otimes_{i+1} W)]
  \ar[dd]|{\eta^{1,i+1}_{X'Z'XZ}\otimes_{j+1} \eta^{1,i+1}_{Y'W'YW}}
  &&[(X'\otimes_{1} X)\otimes_{j+1} (Y'\otimes_{1} Y)]\otimes_{i+1} [(Z'\otimes_{1} Z)\otimes_{j+1} (W'\otimes_{1} W)]
  \ar[dd]|{[M_{AA'A''}\otimes_{j+1} M_{BB'B''}]\otimes_{i+1} [M_{CC'C''}\otimes_{j+1} M_{DD'D''}]}
  \ar[ddll]|{\eta^{i+1,j+1}_{(X'\otimes_{1} X)(Y'\otimes_{1} Y)(Z'\otimes_{1} Z)(W'\otimes_{1} W)}}
  \\\\
  [(X'\otimes_{1} X)\otimes_{i+1} (Z'\otimes_{1} Z)]\otimes_{j+1} [(Y'\otimes_{1} Y)\otimes_{i+1} (W'\otimes_{1} W)]
  \ar[ddr]|{[M_{AA'A''}\otimes_{i+1} M_{CC'C''}]\otimes_{j+1} [M_{BB'B''}\otimes_{i+1} M_{DD'D''}]}
  &&[{\cal A}(A,A'')\otimes_{j+1} {\cal B}(B,B'')]\otimes_{i+1} [{\cal C}(C,C'')\otimes_{j+1} {\cal D}(D,D'')]
  \ar[ddl]|{\eta^{i+1,j+1}_{{\cal A}(A,A''){\cal B}(B,B''){\cal C}(C,C''){\cal D}(D,D'')}}\\\\
  &[{\cal A}(A,A'')\otimes_{i+1} {\cal C}(C,C'')]\otimes_{j+1} [{\cal B}(B,B'')\otimes_{i+1} {\cal D}(D,D'')]
  }
  $$
  \end{small}
  \end{sideways}
  }
  \end{center}\vspace{-3mm}
\nocolon
\caption{}\label{diag2}
\end{figure}

  The lower quadrilateral in Figure \ref{diag2} commutes by naturality
  of $\eta$ and the upper hexagon commutes since it is an instance of
  the giant hexagonal interchange.  As for $\alpha^{(1)}$, the
  compatibility with the unit of $\eta^{(1)i,j}$ follows directly from
  the naturality of $\eta^{i+1,j+1}$ and the fact that
  $j_{[(A,B),(C,D)]} = [(j_A\otimes_{j+1} j_B)\otimes_{i+1}
  (j_C\otimes_{j+1} j_D)]$.

  Also the 2--naturality of $\eta^{(1)i,j}$ follows directly from the
  naturality of $\eta^{i+1,j+1}$ and the Mac Lane coherence theorem.
  
  Since $\alpha^{(1)}$ and $\eta^{(1)}$ are both defined based upon
    $\alpha$ and $\eta$ their ${\cal V}$--functor components satisfy
    all the axioms of the definition of a $k$--fold monoidal category.
    At this level of course it is actually a $k$--fold monoidal
    2--category.
  
  Notice that we have used all the axioms of a $k$--fold monoidal category. The external and internal unit conditions
  imply the unital nature of ${\cal V}$--Cat and the unit axioms for a product of ${\cal V}$--categories respectively.
  The external and internal associativities give us respectively the ${\cal V}$--functoriality of $\alpha^{(1)}$ and 
  the associativity of the composition morphisms for products of ${\cal V}$--categories. This reflects the dual nature
  of the latter two axioms that was pointed out for the braided case in \cite{forcey1}. Finally the giant hexagon gives us precisely the 
  ${\cal V}$--functoriality of $\eta^{(1)}.$ Notice also that we have used in each case the instance of the axiom 
  corresponding to $i=1; j=2\dots k.$  The remaining instances will be used as we iterate the categorical delooping. 
   \end{proof}

\section{Further questions}

For ${\cal V}$ $k$--fold monoidal we have demonstrated that ${\cal V}$--Cat is $(k-1)$--fold monoidal.
    By induction we have that this process continues, i.e. that ${\cal V}$--$n$--Cat $= {\cal V}$--$(n-1)$--Cat--Cat is $(k-n)$--fold
 monoidal for $k > n$. For example, let us expand our description of the next level: the fact that ${\cal V}$--2--Cat = ${\cal V}$--Cat--Cat 
is $(k-2)$--fold monoidal. Now we are considering enrichment over ${\cal V}$--Cat. All the constructions in the proof above are recursively repeated.
The unit ${\cal V}$--2--category is  
  denoted as ${\bcal I}$ where $\left|{\bcal I}\right| = \{\mathbf{0}\}$ and ${\bcal I}(\mathbf{0},\mathbf{0}) = {\cal I}.$
  Products of ${\cal V}$--2--categories are given by ${\bcal U} \otimes^{(2)}_i {\bcal W}$ for $i=1\ldots k-2.$ Objects are pairs of objects as usual,
  and that there are exactly $k-2$ products is seen when the definition of hom--objects is given. In ${\cal V}$--2--Cat,
  $$ [{\bcal U} \otimes^{(2)}_i {\bcal W}]((U,W),(U',W')) = {\bcal U}(U,U') \otimes^{(1)}_{i+1} {\bcal W}(W,W'). $$
  Thus we have that
  $$ [{\bcal U} \otimes^{(2)}_i {\bcal W}]((U,W),(U',W'))((f,f'),(g,g')) $$
  $$ = [{\bcal U}(U,U') \otimes^{(1)}_{i+1} {\bcal W}(W,W')]((f,f'),(g,g')) $$
  $$ = {\bcal U}(U,U')(f,g) \otimes_{i+2} {\bcal W}(W,W')(f',g'). $$
  The definitions of $\alpha^{(2)i}$ and $\eta^{(2)i,j}$  are just as in the lower
  case. 
  For instance, $\alpha^{(2)i}$ will now be a 3--natural transformation, that is, a family of
  ${\cal V}$--2--functors 
  $$\alpha^{(2)i}_{{\bbcal U}{\bbcal V}{\bbcal W}}:({\bcal
  U}\otimes^{(2)}_i {\bcal V}) \otimes^{(2)}_i {\bcal W} \to {\bcal
  U}\otimes^{(2)}_i ({\bcal V}\otimes^{(2)}_i {\bcal W}).$$ To each of
  these is associated a family of ${\cal V}$--functors
  $$\alpha^{(2)i}_{{\bbcal U}{\bbcal V}{\bbcal W}_{(U,V,W)(U',V',W')}} = \alpha^{(1)i+1}_{{\bbcal U}(U,U'){\bbcal V}(V,V'){\bbcal W}(W,W')}$$
  to each of which is associated a family of hom--object morphisms:
  $$\alpha^{(2)i}_{{\bbcal U}{\bbcal V}{\bbcal W}_{(U,V,W)(U',V',W')_{(f,g,h)(f',g',h')}}} = \alpha^{i+2}_{{\bbcal U}(U,U')(f,f'){\bbcal V}(V,V')(g,g'){\bbcal W}(W,W')(h,h')}$$
  Verifications that these define
  a valid $(k-2)$--fold monoidal 3--category all follow just as in the lower dimensional
  case. 
  The facts about the ${\cal V}$--functors are shown by using the original $k$--fold monoidal category axioms that involve i=2.

 In the next paper \cite{Forcey} my aim is to show how enrichment increases categorical dimension as it 
decreases monoidalness. That paper also includes the definitions of ${\cal V}$--$n$--categories and of the
morphisms of ${\cal V}$--$n$--Cat. In further work I want to relate
enrichment more precisely to topological delooping as well as to other categorical constructions that have similar
topological implications. 

In \cite{StAlg} Street defines the nerve of a strict $n$--category.
    Recently Duskin in \cite{Dusk} has worked out the description of the nerve of a bicategory. 
    This allows us to ask whether
    these nerves will prove to be the logical link to loop spaces for higher dimensional iterated monoidal categories. 
        
Passing to the category of enriched categories basically reduces the number
    of products so that for ${\cal V}$ a $k$--fold monoidal $n$--category, ${\cal V}$--Cat becomes a $(k-1)$--fold 
    monoidal $(n+1)$--category. This picture was anticipated by Baez and Dolan \cite{Baez1} in the context where the 
    $k$--fold monoidal $n$--category is specifically a (weak) $(n+k)$--category with only one object,
    one 1--cell, etc. up to only one $k$--cell. Their version of categorical delooping simply consists of creating  
    from a monoidal category ${\cal V}$ the one object bicategory that has its morphisms the objects of ${\cal V}.$
    Relating the two versions of delooping is important to an understanding of how categories model spaces.

\Addresses\recd
}}}
\end{document}